\documentclass[12pt]{article}

\usepackage{amsfonts}

\newtheorem{thm}{Theorem}[section]
\newtheorem{lemma}[thm]{Lemma}
\newtheorem{cor}[thm]{Corollary}

\newtheorem{prop}[thm]{Proposition}
\newtheorem{conjecture}{Conjecture}

\newcommand{\proof
}{\par\medskip\noindent {\bf Proof.\ \ }}

\newcommand{\be}{\begin{equation}}
\newcommand{\ee}{\end{equation}}
\newcommand{\openbox}{\leavevmode
  \hbox to8pt{\hfil\vrule\vbox to6pt{\hrule width6pt\vfil\hrule}\vrule}}

\newcommand{\qed}{\hbox to5pt{ } \hfill \openbox\bigskip\medskip}

\newcommand{\cT}{\mbox{$\cal T$}}

\newcommand{\cK}{\mbox{$\cal K$}}
\newcommand{\cF}{\mbox{$\cal F$}}
\newcommand{\cG}{\mbox{$\cal G$}}
\newcommand{\cH}{\mbox{$\cal H$}}

\newcommand{\cM}{\mbox{$\cal M$}}

\newcommand{\Z}{\mathbb Z}

\title{About  sunflowers}
\author{G\'abor Heged\H{u}s
\\{\normalsize  \'Obuda University}
\\{\normalsize B\'ecsi \'ut 96, Budapest, Hungary, H-1037}
\\{\normalsize hegedus.gabor@nik.uni-obuda.hu}
}

\begin{document}
\maketitle

\begin{abstract}
Alon,  Shpilka and Umans considered the following version of usual sunflower-free subset:  a subset $\mbox{$\cal F$}\subseteq  \{1,\ldots ,D\}^n$ for $D>2$ is sunflower-free if for every distinct triple  $x,y,z\in \mbox{$\cal F$}$ there exists a coordinate $i$ where exactly two of $x_i,y_i,z_i$ are equal. Combining the polynomial method with character theory Naslund and Sawin proved that any sunflower-free set $\mbox{$\cal F$}\subseteq  \{1,\ldots ,D\}^n$ has size
$$
|\mbox{$\cal F$}|\leq c_D^n,
$$
where $c_D=\frac{3}{2^{2/3}}(D-1)^{2/3}$. 

In this short note  we  give a new upper bound for the size of  sunflower-free subsets of  $\{1,\ldots ,D\}^n$.

Our main result is a new upper bound for the size of  sunflower-free $k$-uniform subsets.
 
More precisely, let $k$ be an arbitrary integer. Let $\mbox{$\cal F$}$ be a sunflower-free $k$-uniform set system. 
Consider
$M:=|\bigcup\limits_{F\in \cF} F|.
$
Then 
$$
|\mbox{$\cal F$}|\leq 3(\lceil\frac{2k}{3}\rceil+1)(2^{1/3}\cdot 3e)^k(\lceil\frac Mk\rceil -1)^{\lceil\frac{2k}{3}\rceil}.
$$

In the proof we use  Naslund and Sawin's result about sunflower-free subsets in $\{1,\ldots ,D\}^n$.
\end{abstract}
\medskip

\noindent

\section{Introduction}

Let  $[n]$ stand for the set $\{1,2,
\ldots, n\}$. We denote the family of all subsets of $[n]$  by $2^{[n]}$. 

Let $X$ be a fixed subset of $[n]$. Let $0\leq k\leq n$ integers. We denote by
${X \choose k}$ the family of all  $k$ element subsets of $X$.

We say that a family $\cF$ of subsets of $[n]$ {\em $k$-uniform}, if $|F|=k$ for each $F\in \cF$.

Recall that a family $\cF=\{F_1,\ldots ,F_m\}$ of subsets of $[n]$ is a {\em sunflower} (or {\em $\Delta$-system}) with $t$ petals if
$$
F_i\cap F_j=\bigcap\limits_{s=1}^t F_s
$$
for each $1\leq i,j\leq t$.

The {\em kernel} of a sunflower is  the intersection of the members of this sunflower.

By definition a family of disjoint sets is a sunflower with empty kernel.

Erd\H{o}s and  Rado gave a remarkable upper bound for the size of a $k$-uniform family without a sunflower with $t$ petals (see \cite{ER}).
\begin{thm} \label{Sthm} (Sunflower theorem)
If $\cF$ is a $k$-uniform set system with more than 
$$
k!(t-1)^k\big( 1-\sum_{s=1}^{k-1} \frac{s}{(s+1)! (t-1)^s}\big)
$$ 
members, then $\cF$ contains a sunflower with $t$ petals.
\end{thm}

Later Kostochka improved this upper bound in \cite{K2}. 
\begin{thm} \label{Kost}
Let $t>2$ and $\alpha>1$ be fixed integers. Let $k$ be an arbitrary integer. Then there exists a constant $D(t,\alpha)$ such that if  $\cF$ is a $k$-uniform set system with more than 
$$
D(t,\alpha)k! \Big( \frac{(\log \log \log k)^2}{\alpha \log \log k}\Big)^k
$$
 members, then $\cF$ contains a sunflower with $t$ petals.
\end{thm}

The following statement is conjectured by Erd\H{o}s and  Rado  in \cite{ER}.
\begin{conjecture} \label{Econj}
For each $t$, there exists a constant $C(t)$ such that  if  $\cF$ is a $k$-uniform set system with more than 
$$
C(t)^k
$$
members,  then $\cF$ contains a sunflower with $t$ petals.
\end{conjecture}
It is well-known that Erd\H{o}s  offered 1000 dollars for the proof or disproof of this conjecture for $t=3$ (see \cite{E}).





Naslund and Sawin proved the following upper bound for the size of a sunflower-free family in \cite{NS}. Their argument based on Tao's slice--rank bounding method (see the blog \cite{T}). 

\begin{thm} \label{NS}
Let $\cF$ be a family of subsets of  $[n]$ without a sunflower with 3 petals. Then
$$
|\cF|\leq 3(n+1)\sum_{i=0}^{\lfloor n/3\rfloor} {n\choose i}.
$$ 
\end{thm}

Alon,  Shpilka and Umans considered the following version of usual sunflowers in \cite{ASU}: 
Let $D>2$, $n\geq 1$ be  integers. Then $k$ vectors $v_1,\ldots , v_k\in {\Z}_D^n$ form a {\em $k$-sunflower} if for every coordinate $i\in [n]$ it holds that either $(v_1)_i=\ldots =(v_k)_i$ or they all differ on that coordinate.   

In the following the 'sunflower' term means always a $3$-sunflower.

Naslund and Sawin gave the following upper bounds for  the size of sunflower-free families in \cite{NS} Theorem 2. Their proof worked only for $3$-sunflowers.
\begin{thm} \label{NS2}
Let $D>2$, $n\geq 1$ be  integers. Let $\cF\subseteq  {\Z}_D^n$ be a sunflower-free family in  ${\Z}_D^n$. Then 
$$
|\cF| \leq c_D^n,
$$
where $c_D=\frac{3}{2^{2/3}}(D-1)^{2/3}$. 
\end{thm}

Let $D>2$, $n\geq 1$ be integers. Let $s(D,n)$ denote the maximum size of a sunflower-free family in  ${\Z}_D^n$. 

Define $J(q):=\frac 1q \Big( \min_{0<x<1} \frac{1-x^q}{1-x} x^{-\frac{q-1}{3}}\Big)$ for each $q>1$.

This $J(q)$  constant  appeared in Ellenberg and Gijswijt's bound for the size of three-term progression-free sets (see \cite{EG}). Blasiak, Church, Cohn, Grochow  and Umans proved in \cite{BCCGU} Proposition 4.12 that $J(q)$ is a decreasing function of $q$ and
$$
\lim_{q\to \infty} J(q)=\inf_{z>3} \frac{z-z^{-2}}{3\log(z)}=0.8414\ldots .
$$ 
It is easy to verify that $J(3)=0.9184$, consequently $J(q)$ lies in the range 
$$
0.8414\leq J(q)\leq 0.9184
$$
for each $q\geq 3$.

Since a sunflower-free family in  ${\Z}_D^n$ can not contain a a three-term arithmetic progression, hence the  Ellenberg and Gijswijt's striking result
(see \cite{EG}) implies the following upper bound.

\begin{thm} \label{EG_up}
Let $n\geq 1$ be an integer,  $p^{\alpha}>2$ be a prime power. Let $\cF\subseteq  {\Z}_D^n$ be a sunflower-free family in  ${\Z}_D^n$. Then 
$$
|\cF| \leq (J(p^{\alpha})p^{\alpha})^n.
$$
\end{thm}

Now we give some new bounds for the size of sunflower-free families in  ${\Z}_D^n$.

The Chinese Remainder Theorem implies immediately the following result.
\begin{thm} \label{EG_up_CRT}
Let $m=p_1^{\alpha_1}\cdot \ldots \cdot p_r^{\alpha_r}$, where $p_i$ are different primes. Then
$$
s(m,n)\leq s(p_1^{\alpha_1},n)\cdot \ldots \cdot s(p_r^{\alpha_r},n)\leq ((\prod_{i=1}^r J(p_i^{\alpha_i}))m)^n.
$$
\end{thm}

\proof

By the Chinese Remainder Theorem there exists a bijection
$$
\phi: {\Z}_m\to {\Z}_{p_1^{\alpha_1}} \times \ldots \times {\Z}_{p_r^{\alpha_r}}.
$$
We can extend this bijection in a natural way to $({\Z}_m)^n$ and we get the bijection
$$
\phi^{*}: ({\Z}_m)^n\to ({\Z}_{p_1^{\alpha_1}})^n \times \ldots \times ({\Z}_{p_r^{\alpha_r}})^n.
$$

Let $\cF\subseteq  ({\Z}_m)^n$ be a sunflower-free family in  $({\Z}_m)^n$. 

Then it is easy to check that $\phi^{*}(\cF)$ is a a sunflower-free family in $({\Z}_{p_1})^n \times \ldots \times ({\Z}_{p_r})^n$. Hence
$$
|\cF|\leq |\phi^{*}(\cF)| \leq   s(p_1^{\alpha_1},n)\cdot \ldots \cdot s(p_r^{\alpha_r},n).
$$
\qed

Next we give an other new upper bound for the size of sunflower-free families in  ${\Z}_D^n$, which is independent of $D$.
\begin{thm} \label{Sf_up_gen}
Let $D>2$, $n\geq 1$ be  integers, $\alpha>1$ be a real number. Let $\cF\subseteq  {\Z}_D^n$ be a sunflower-free family in  ${\Z}_D^n$. Then there exists a constant $K(\alpha)>0$ such that
$$
|\cF|\leq K(\alpha) n! \Big( \frac{(\log \log \log n)^2}{\alpha \log \log n}\Big)^n.
$$
\end{thm}

Our main result is a new upper bound for the size of $k$-uniform  sunflower-free families. In the proof we use  Theorem \ref{NS2} and Erd\H{o}s and Kleitman's famous result about $k$-partite hypergraphs.

\begin{thm}  \label{USf_up_gen2}
Let $k$ be an arbitrary integer. Let $\cF$ be a sunflower-free $k$-uniform set system. Let 
$M:=|\bigcup\limits_{F\in \cF} F|.
$
Then 
$$
|\cF|\leq 3(\lceil\frac{2k}{3}\rceil+1)(2^{1/3}\cdot 3e)^k(\lceil\frac Mk\rceil -1)^{\lceil\frac{2k}{3}\rceil}.
$$
\end{thm}




\begin{cor} \label{USf_up_gen}
Let $k$ be an arbitrary integer. Let $\epsilon>0$ be a fixed real number. Let $\cF$ be a sunflower-free $k$-uniform set system. Suppose that 
$$
|\bigcup\limits_{F\in \cF} F|\leq k^{2.5-\epsilon}.
$$
Then 
$$
|\cF|\leq 3(\lceil\frac{2k}{3}\rceil+1)(2^{1/3}\cdot 3e)^k k^{\lceil k(1-\frac{2\epsilon}{3})\rceil}.
$$
\end{cor}
\proof

Define 
$M:=|\bigcup\limits_{F\in \cF} F|$. Then
$$
\frac Mk\leq k^{1.5-\epsilon}.
$$

Theorem \ref{USf_up_gen2} gives us the desired result.
\qed

We present our proofs in Section 2.

\section{Proofs}

{\bf Proof of Theorem \ref{Sf_up_gen}:}\\

Let $\cF\subseteq  ({\Z}_D)^n$ be a sunflower-free family in  $({\Z}_D)^n$. 
We define first a hypergraph corresponding to $\cF$.

Let $U:=[D]\times [n]$ denote the universe of this hypergraph.

Then for each vector $v\in {\Z}_D^n$ we can define the set 
$$
M(v):=\{(v_1+1,1), \ldots ,  (v_n+1,n)\}\subseteq U.
$$

It is clear that $M(v)$ are $n$-sets.
Consider the hypergraph
$$
M(\cF):=\{M(v):~ v\in \cF\}.
$$

Then $M(\cF)$ is an $n$-uniform set family. 

It is easy to check that $M(\cF)$ is a sunflower-free hypergraph, since $\cF$ is sunflower-free. Consequently we can apply Theorem \ref{Kost} to the hypergraph $M(\cF)$ and we get our result.

\qed

Suppose that $\cK\subseteq {X \choose k}$ and that for some disjoint decomposition 
$$
X=X_1\oplus \ldots \oplus X_m,
$$
$\cK$ satisfies the equality $|F\cap X_i|=1$ for all $F\in \cK$ and  $1\leq i\leq m$. Then $\cK$ is an $m$-partite hypergraph.

Erd\H{o}s and Kleitman proved in \cite{EK} the following well-known result using  an averaging argument.

\begin{thm} \label{EK}
 Suppose that $\cK\subseteq {X \choose k}$. Then there exists a subfamily $\cG\subseteq \cF$ such that 
$\cG$ is $k$-partite and satisfies 
\begin{equation} \label{upper_EK}
|\cG|\geq \frac{k!}{k^k} |\cK|.
\end{equation}
\end{thm}

We use also in our proof the following generalization of Theorem
\ref{NS2}.

\begin{thm} \label{NS_gen} 
Let $D_i\geq 3$ be integers for each $i\in [n]$. 
Let $\cH\subseteq {\Z}_{D_1} \times \ldots \times {\Z}_{D_n}$ be a sunflower-free family in ${\Z}_{D_1} \times \ldots \times {\Z}_{D_n}$. Then 
$$
|\cH|\leq 3\sum_{I\subseteq [n], 0\leq |I|\leq \frac{2n}{3}} \prod_{i\in I} (D_i-1).
$$
\end{thm}

\proof

A simple modification of the argument appearing  the proof of Theorem \ref{NS2} works as a proof of  Theorem  \ref{NS_gen}. 
\qed

It is easy to verify the following Proposition.

\begin{prop} \label{balanced} 
Let $D_i\geq 3$ be integers for each $i\in [n]$. Define $M:=\sum_i D_i$. Then 
$$
\sum_{I\subseteq [n], 0\leq |I|\leq \frac{2n}{3}} \prod_{i\in I} (D_i-1) \leq \sum_{j=0}^{\frac{2n}{3}} {n\choose j}\Big( \lceil\frac Mn \rceil -1 \Big)^j.
$$
\end{prop}
\qed

{\bf Proof of Theorem \ref{USf_up_gen}:}\\

Let $\cF$ be a sunflower-free $k$-uniform set system. 

Define $M:= |\bigcup\limits_{F\in \cF} F|$ and $X:=\bigcup\limits_{F\in \cF} F$.

By  Theorem  \ref{EK} there exists a subfamily $\cG \subseteq \cF$ such that 
$\cG$ is $k$-partite and satisfies 
$$
|\cG|\geq \frac{k!}{k^k} |\cF|\geq  \frac{|\cF|}{e^{k}}.
$$

Consider the disjoint decomposition  into classes 
$$
X=C_1\oplus \ldots \oplus C_k,
$$
where $\cG$ satisfies the equality $|G\cap C_i|=1$ for all $G\in \cG$ and  $1\leq i\leq k$.

We can suppose that $C_1, \ldots ,C_t$ are the classes with $|C_i|= 2$ for each  $1\leq i\leq t$ and   $|C_i|\geq 3$  for each $i>t$.

Let $C_i=\{x_i, y_i\}$ for each $1\leq i\leq t$.  Define for each $N\subseteq [t]$ the following subfamily of $\cG$:
$$
\cG(N):=\{G\in \cG:~ \{x_i:~ i\in N\}\cup \{y_i:~ i\in [t]\setminus N\}\subseteq G\}.
$$

Denote by $L\subseteq [t]$ the subset with
$$
|\cG(L)|=\max_{N\subseteq [t]} |\cG(N)|.
$$
Consider the set system
$$
\cH:=\{F\setminus L:F\in \cG(L)\}.
$$

Clearly here $X\setminus L=\bigcup\limits_{H\in \cH} H$.

Then $\cH$ is a $(k-t)$-uniform, $(k-t)$-partite set system with the disjoint decomposition  into classes
$$
X\setminus L=B_1\oplus \ldots \oplus B_{k-t},
$$
where $\cH$ satisfies the equality $|H\cap B_i|=1$ for all $H\in \cH$ and  $1\leq i\leq k-t$. Our construction of the set system $\cH$ shows that $|B_i|\geq 3$ for each $1\leq i\leq k-t$.

On the other hand it follows from the equality 
$$
|\cG(L)|=\max_{N\subseteq [t]} |\cG(N)|
$$
that
\begin{equation} \label{upper2}
|\cG|\leq \sum_{N\subseteq [t]} |\cG(N)| \leq 2^t |\cG(L)|=2^t  |\cH|\leq 2^k |\cH|.
\end{equation}

In the following we consider only the case when  $t=0$. The $t>0$ cases can be treated in a similar way.

We use  the following Proposition in our proof.
\begin{prop} \label{bij} 
Let $D_i\geq 3$ be integers for each $i\in [n]$. Define $M:=\sum_i D_i$. Then there exists an injection $\psi:{\Z}_{D_1} \times \ldots \times {\Z}_{D_n}\to {[M]\choose n}$ such that each $n$-uniform, $n$-partite  set system with classes $C_i$, where $|C_i|=D_i$ for each $1\leq i\leq n$ is precisely the image set of the map $\psi$ and each  $n$-uniform, $n$-partite and  sunflower-free   family with classes $C_i$, where $|C_i|=D_i$ for each $1\leq i\leq n$ corresponds to a sunflower-free family in  ${\Z}_{D_1} \times \ldots \times {\Z}_{D_n}$.
\end{prop}

We can apply Proposition \ref{bij} with the choices $D_i:=|B_i|$ and we get that $\cT:=\psi^{-1}(\cH)$ is a sunflower-free family in  the group ${\Z}_{D_1} \times \ldots \times {\Z}_{D_k}$.

It follows from Theorem \ref{NS_gen}  that
$$
|\cH|=|\cT|\leq 3\cdot \sum_{I\subseteq [n], 0\leq |I|\leq \frac{2k}{3}} \prod_{i\in I} (D_i-1).
$$

But we get from Proposition \ref{balanced}  that 
$$
3\cdot\sum_{I\subseteq [n], 0\leq |I|\leq \frac{2k}{3}} \prod_{i\in I} (D_i-1) \leq 3\sum_{j=0}^{\frac{2k}{3}} {k\choose j}\Big( \lceil\frac  Mk\rceil -1 \Big)^j.
$$

Hence
$$
3\sum_{j=0}^{\frac{2k}{3}} {k\choose j}\Big(\lceil \frac Mk \rceil-1 \Big)^j\leq 3(\lceil\frac{2k}{3}\rceil+1)(\frac{3}{2^{2/3}})^k(\frac Mk -1)^{\lceil\frac{2k}{3}\rceil}.
$$

The desired upper bound follows from equations (\ref{upper_EK}) and (\ref{upper2}):
$$
|\cF|\leq e^k|\cG|\leq 2^ke^k|\cH|\leq 3(\lceil\frac{2k}{3}\rceil+1)(2^{1/3}\cdot 3e)^k(\lceil\frac Mk\rceil -1)^{\lceil\frac{2k}{3}\rceil}.
$$
\qed

\section{Concluding remarks}

The following conjecture implies an unconditional, strong upper bound for the size of any sunflower-free $k$-uniform set system.
\begin{conjecture} \label{Hconj}
There exists a $D>0$ constant such that if $\cF$ is any sunflower-free $k$-uniform set system, then 
$$
|\bigcup\limits_{F\in \cF} F|\leq Dk^2.
$$
\end{conjecture}

We give here a weaker version of Conjecture \ref{Hconj}.
\begin{conjecture} \label{Hconj2}
Let $\cF$ be a sunflower-free $k$-uniform set system. Then there exist $F_1,\ldots ,F_{2k}\in \cF$ such that
$$
\bigcup\limits_{F\in \cF} F=\bigcup\limits_{j=1}^{2k} F_j.
$$
\end{conjecture}



\end{document}